\documentclass[11pt]{article}
\usepackage[centertags]{amsmath}
\usepackage{amsfonts}
\usepackage{amssymb}
\usepackage{amsthm}
\usepackage{newlfont}
\usepackage{graphicx}

\newfont{\bb}{msbm10 at 12pt}
\def\r{\hbox{\bb R}}

\def\e{\hbox{\bf E}}
\def\t{\hbox{\bf T}}
\def\n{\hbox{\bf N}}
\def\b{\hbox{\bf B}}

\newtheorem{theorem}{Theorem}[section]
\newtheorem{definition}[theorem]{Definition}

\newtheorem{lemma}[theorem]{Lemma}

\begin{document}
\title{Similar Curves With Variable Transformations}

\author{Mostafa F. El-Sabbagh\\
Mathematics Department, Faculty of Science,\\
Minia University, Minia,  Egypt.\\
E-mail: \textit{msabbagh52@yahoo.com}\\
\centerline{}\\
Ahmad T. Ali \footnote{Corresponding author.}\\
Mathematics Department,\\
Faculty of Science, Al-Azhar University,\\
Nasr City, 11448, Cairo, Egypt.\\
E-mail: \textit{atali71@yahoo.com}}

\maketitle
\begin{abstract}
In this paper, we define a new family of curves and call it a {\it family of similar curves with variable transformation} or briefly {\it SA-curves}. Also we introduce some characterizations of this family and we give some theorems. This definition introduces a new classification of a space curve. Also, we use this definition to deduce the position vectors of plane curves, general helices and slant helices, as examples of a similar curves with variable transformation.\\

\textbf{Math. Sub. Class. 2000}: 53C40, 53C50\\
\textbf{Keywords}: Classical differential geometry; general helix; slant helices; Intrinsic equations, similar curves.
\end{abstract}

\section{Introduction}
In the local differential geometry, we think of a curve as a geometric set of points, or locus. Intuitively, we are thinking of a curve as the path traced out by a particle moving in $\e^3$. So, investigating position vector of the curve is a classical aim to determine behavior of the particle (curve).

From the view of differential geometry, a {\it straight line} is a geometric curve with the curvature $\kappa(s)=0$. A {\it plane curve} is a family of geometric curves with torsion $\tau(s)=0$. Helix is a geometric curve with non-vanishing constant curvature $\kappa$ and non-vanishing constant torsion $\tau$ \cite{barros}. The helix may be called a {\it circular helix} or {\it $W$-curve} \cite{ilarslan}. It is known that straight line ($\kappa(s)=0$) and circle ($\kappa(s)=a,\,\tau(s)=0$) are degenerate-helices examples \cite{kuhn}. In fact, circular helix is the simplest three-dimensional spirals \cite{camci}.

A curve of constant slope or {\it general helix} in Euclidean 3-space $\e^3$ is defined by the property that the tangent makes a constant angle with a fixed straight line called the axis of the general helix. A classical result stated by Lancret in 1802 and first proved by de Saint Venant in 1845 (see \cite{struik} for details) says that: {\it A necessary and sufficient condition that a curve be a general helix is that the ratio $$\dfrac{\tau}{\kappa}$$ is constant along the curve, where $\kappa$ and $\tau$ denote the curvature and the torsion, respectively}. General helices or {\it inclined curves} are well known curves in classical differential geometry of space curves \cite{milm} and we refer to the reader for recent
works on this type of curves \cite{ali1, gluck, mont2, turgut}.

Izumiya and Takeuchi \cite{izumi} have introduced the concept of {\it slant helix} by saying that the normal lines make a constant angle with a fixed straight line. They characterize a slant helix if and only if the {\it geodesic curvature} of the principal image of the principal normal indicatrix
$$
\sigma=\frac{\kappa^2}{(\kappa^2+\tau^2)^{3/2}}\Big(\frac{\tau}{\kappa}\Big)'
$$
is a constant function. Kula and Yayli \cite{kula1} have studied spherical images of tangent indicatrix and binormal indicatrix of a slant helix and they showed that the spherical images are spherical helices. Recently, Kula et al. \cite{kula2} investigated the relation between a general helix and a slant helix. Moreover, they obtained some differential equations which are characterizations for a space curve to be a slant helix.

A family of curves with constant curvature but non-constant torsion is called Salkowski curves and a family of curves with constant torsion but non-constant curvature is called anti-Salkowski curves $\cite{salkow}$. Monterde \cite{mont1} studied some characterizations of these curves and he proved that the principal normal vector makes a constant angle with fixed straight line. So that: Salkowski and anti-Salkowski curves are the important examples of slant helices.

A unit speed curve of {\it constant precession} in Euclidean 3-space $\e^3$ is defined by the property that its (Frenet) Darboux vector
$$
W=\tau\,\t+\kappa\,\b
$$
revolves about a fixed line in space with constant angle and constant speed. A curve of constant precession is characterized by having
$$
\kappa=\frac{\mu}{m}\sin[\mu\,s],\,\,\,\,\,\,\,\,\,\,\tau=\frac{\mu}{m}\cos[\mu\,s]
$$
or
$$
\kappa=\frac{\mu}{m}\cos[\mu\,s],\,\,\,\,\,\,\,\,\,\,\tau=\frac{\mu}{m}\sin[\mu\,s]
$$
where $\mu$ and $m$ are constants. This curve lies on a circular one-sheeted hyperboloid
$$
x^2+y^2-m^2\,z^2=4m^2
$$
The curve of constant precession is closed if and only if $n=\frac{m}{\sqrt{1+m^2}}$ is rational \cite{scofield}. Kula and Yayli \cite{kula1} proved that the geodesic curvature of the spherical image of the principal normal indicatrix of a curve of constant precession is a constant function equals $-m$. So we may say that: the curves of constant precessions are the important examples of slant helices.

Many important results in the theory of curves in $\e^3$ were initiated by G. Monge and G. Darboux pioneered the moving frame idea. Thereafter, F. Frenet defined his moving frame and his special equations which play important role in mechanics and kinematics as well as in differential geometry \cite{boyer}.

In this work, we define a new family of curves and we call it a family of {\it similar curves with variable transformation} or in brief {\it SA-curves}. Also, we introduce some characterizations of this family and give some theorems. This definition introduces a new classification of a space curve. In the last of this paper, we use this definition to deduce the position vectors of some important special curves. We hope these results will be helpful to mathematicians who are specialized on mathematical modeling as well as other applications of interest.

\section{Preliminaries }
In Euclidean space $\e^3$, it is well known that to each unit speed curve with at least four continuous derivatives, one can associate three mutually orthogonal unit vector fields $\t$, $\n$ and $\b$ are respectively, the tangent, the principal normal and the binormal vector fields \cite{hacis}.

We consider the usual metric in Euclidean 3-space $\e^3$, that is,
$$
\langle,\rangle=dx_1^2+dx_2^2+dx_3^2,
$$
where $(x_1,x_2,x_3)$ is a rectangular coordinate system of $\e^3$.  Let $\psi:I\subset\r\rightarrow\e^3$, $\psi=\psi(s)$, be an arbitrary curve in $\e^3$. The curve $\psi$ is said to be of unit speed (or parameterized by the  arc-length) if $\langle\psi'(s),\psi'(s)\rangle=1$ for any $s\in I$. In particular, if $\psi(s)\not=0$ for any $s$, then it is possible to re-parameterize $\psi$, that is, $\alpha=\psi(\phi(s))$ so that $\alpha$ is parameterized by the arc-length. Thus, we will assume throughout this work that $\psi$ is a unit speed curve.

Let $\{\t(s),\n(s),\b(s)\}$ be the moving frame along $\psi$, where the vectors $\t, \n$ and $\b$ are mutually orthogonal vectors satisfying $\langle\t,\t\rangle=\langle\n,\n\rangle=\langle\b,\b\rangle=1$.
The Frenet equations for $\psi$ are given by (\cite{struik})
\begin{equation}\label{u1}
 \left[
   \begin{array}{c}
     \t'(s) \\
     \n'(s) \\
     \b'(s) \\
   \end{array}
 \right]=\left[
           \begin{array}{ccc}
             0 & \kappa(s) & 0 \\
             -\kappa(s) & 0 & \tau(s) \\
             0 & -\tau(s) & 0 \\
           \end{array}
         \right]\left[
   \begin{array}{c}
     \t(s) \\
     \n(s) \\
     \b(s) \\
   \end{array}
 \right].
 \end{equation}

If $\tau(s)=0$ for all $s\in I$, then $\b(s)$ is a constant vector $V$ and the curve $\psi$ lies in a $2$-dimensional affine subspace orthogonal to $V$, which is isometric to the Euclidean $2$-space $\e^2$.

\section{Position vector of a space curve}
The problem of the determination of parametric representation of the position vector of an arbitrary space curve according to its intrinsic equations is still open in the Euclidean space $\e^3$ \cite{eisenh, lips}. This problem is not easy to solve in general case. However, this problem is solved in three special cases only, Firstly, in the case of a plane curve ($\tau=0$). Secondly, in the case of a helix ($\kappa$ and $\tau$ are both non-vanishing constant). Recently, Ali \cite{ali2, ali3} adapted fundamental existence and uniqueness theorem for space curves in Euclidean space $\e^3$ and constructed a vector differential equation to solve this problem in the case of a general helix ($\frac{\tau}{\kappa}$ is constant) and in the case of a slant helix
\begin{equation}\label{u02}
\frac{\tau(s)}{\kappa(s)}=\pm\dfrac{m\,\int\kappa(s)ds}{\sqrt{1-m^2\Big(\int\kappa(s)ds\Big)^2}},
\end{equation}
where $m=\frac{n}{\sqrt{1-n^2}}$, $n=\cos[\phi]$ and $\phi$ is the constant angle between the axis of a slant helix and the principal normal vector. However, this problem is not solved in other cases of space curves.\\

Now we describe this problem within the following theorem:

\begin{theorem}\label{th01} Let $\psi=\psi(s)$ be an unit speed curve parameterized by the arclength $s$. Suppose $\psi=\psi(\theta)$ is another parametric representation of this curve by the parameter $\theta=\int\kappa(s)ds$.  Then, the tangent vector $\t$ satisfies a vector differential equation of third order as follows:
\begin{equation}\label{u03}
\Big(\frac{1}{f(\theta)}\t''(\theta)\Big)'+\Big(\frac{1+f^2(\theta)}{f(\theta)}\Big)\t'(\theta)
-\frac{f'(\theta)}{f(\theta)}\t(\theta)=0,
\end{equation}
where $f(\theta)=\frac{\tau(\theta)}{\kappa(\theta)}$.
\end{theorem}

{\bf Proof.} Let $\psi=\psi(s)$ be a unit speed curve. If we write this curve in another parametric representation $\psi=\psi(\theta)$, where $\theta=\int\kappa(s)ds$, we have new Frenet equations as follows:
\begin{equation}\label{u04}
 \left[
   \begin{array}{c}
     \t'(\theta) \\
     \n'(\theta) \\
     \b'(\theta) \\
   \end{array}
 \right]=\left[
           \begin{array}{ccc}
             0 & 1 & 0 \\
             -1 & 0 & f(\theta) \\
             0 & -f(\theta) & 0 \\
           \end{array}
         \right]\left[
   \begin{array}{c}
     \t(\theta) \\
     \n(\theta) \\
     \b(\theta) \\
   \end{array}
 \right],
 \end{equation}
where $f(\theta)=\frac{\tau(\theta)}{\kappa(\theta)}$. If we substitute the first equation of the new Frenet equations (\ref{u04}) to second equation of (\ref{u04}), we have
\begin{equation}\label{u05}
\b(\theta)=\frac{1}{f(\theta)}\Big[\t''(\theta)+\t(\theta)\Big].
\end{equation}
Substituting the above equation in the last equation from (\ref{u04}), we obtain a vector differential equation of third order (\ref{u03}) as desired.

The equation (\ref{u03}) is not easy to solve in general case. If one solves this equation, the natural representation of the position vector of an arbitrary space curve can be determined as follows:
\begin{equation}\label{u06}
\psi(s)=\int\,\t(s)\,ds+C,
\end{equation}
or in parametric representation
\begin{equation}\label{u07}
\psi(\theta)=\int\,\frac{1}{\kappa(\theta)}\,\t(\theta)\,d\theta+C,
\end{equation}
where $\theta=\int\kappa(s)ds$ and $C$ is a constant vector.

\section{Similar curves with variable transformations}

\begin{definition}\label{df-main1} Let $\psi_\alpha(s_\alpha)$ and $\psi_\beta(s_\beta)$ be two regular curves in $\e^3$ parameterized by arclengths $s_\alpha$ and $s_\beta$ with curvatures $\kappa_\alpha$ and $\kappa_\beta$, torsion $\tau_\alpha$ and $\tau_\beta$ and Frenet frames $\{\t_\alpha,\n_\alpha,\b_\alpha\}$ and $\{\t_\beta,\n_\beta,\b_\beta\}$. $\psi_\alpha(s_\alpha)$ and $\psi_\beta(s_\beta)$ are called similar curves with variable transformation $\lambda^\alpha_\beta$ if there exist a variable transformation
\begin{equation}\label{u08}
s_\alpha=\int\,\lambda^\alpha_\beta(s_\beta)\,ds_\beta
\end{equation}
of the arclengths such that the tangent vectors are the same for the two curves i.e.,
\begin{equation}\label{u3}
\t_\beta(s_\beta)=\t_\alpha(s_\alpha),
\end{equation}
\end{definition}
for all corresponding values of parameters under the transformation $\lambda^\alpha_\beta$.

Here $\lambda^\alpha_\beta$ is arbitrary function of the arclength. It is worth noting that $\lambda^\alpha_\beta\lambda^\beta_\alpha=1$. All curves satisfy equation (\ref{u3}) are called a family of similar curves with variable transformations. If we integrate the equation (\ref{u3}) we have the following important theorem:

\begin{theorem}\label{th02} The position vectors of the family of similar curves with variable transformations can be written in the following form:
\begin{equation}\label{u010}
\psi_\beta(s_\beta)=\int\,\t_\alpha\Big(s_\alpha(s_\beta)\Big)\,ds_\beta=
\int\,\t_\alpha\Big(s_\alpha\Big)\,\lambda^\beta_\alpha\,ds_\beta.
\end{equation}
\end{theorem}

\begin{theorem}\label{th-03} Let $\psi_\alpha(s_\alpha)$ and $\psi_\beta(s_\beta)$ be two regular curves in $\e^3$. Then $\psi_\alpha(s_\alpha)$ and $\psi_\beta(s_\beta)$ are similar curves with variable transformation if and only if the principal normal vectors are the same for all curves
\begin{equation}\label{u4}
\n_\beta(s_\beta)=\n_\alpha(s_\alpha),
\end{equation}
under the particular variable transformation
\begin{equation}\label{u5}
\lambda^\beta_\alpha=\frac{ds_\beta}{ds_\alpha}=\frac{\kappa_\alpha}{\kappa_\beta}
\end{equation}
of the arc-lengths.
\end{theorem}

{\bf Proof.} ($\Rightarrow$) Let $\psi_\alpha(s_\alpha)$ and $\psi_\beta(s_\beta)$ be two regular similar curves with variable transformations in $\e^3$. Differentiating the equation (\ref{u3}) with respect to $s_\beta$ we have
\begin{equation}\label{u6}
\kappa_\beta(s_\beta)\n_\beta(s_\beta)=\kappa_\alpha(s_\alpha)\n_\alpha(s_\alpha)\frac{ds_\alpha}{ds_\beta}.
\end{equation}
The above equation leads to the two equations (\ref{u4}) and (\ref{u5}).

($\Leftarrow$) Let $\psi_\alpha(s_\alpha)$ and $\psi_\beta(s_\beta)$ be two regular curves in $\e^3$ satisfying the two equations (\ref{u4}) and (\ref{u5}). If we multiply equation (\ref{u5}) by $\kappa_\beta(s_\beta)$ and integrate the result with respect to $s_\beta$ we have
\begin{equation}\label{u8}
\int\,\kappa_\beta(s_\beta)\,\n_\beta(s_\beta)\,ds_\beta=\int\,\kappa_\beta(s_\beta)\,\n_\beta(s_\beta)\,\frac{ds_\beta}{ds_\alpha}\,ds_\alpha.
\end{equation}
From the equations (\ref{u4}) and (\ref{u5}), equation (\ref{u8}) takes the form
\begin{equation}\label{u10}
\int\,\kappa_\beta(s_\beta)\,\n_\beta(s_\beta)\,ds_\beta=\int\,\kappa_\alpha(s_\alpha)\,\n_\alpha(s_\alpha)\,ds_\alpha,
\end{equation}
which leads to (\ref{u3}) and the proof is complete.

\begin{theorem}\label{th04} Let $\psi_\alpha(s_\alpha)$ and $\psi_\beta(s_\beta)$ be two regular curves in $\e^3$. Then $\psi_\alpha(s_\alpha)$ and $\psi_\beta(s_\beta)$ are similar curves with variable transformation if and only if the binormal vectors are the same, i.e.,
\begin{equation}\label{u11}
\b_\beta(s_\beta)=\b_\alpha(s_\alpha),
\end{equation}
\end{theorem}
under arbitrary variable transformations $s_\beta=s_\beta(s_\alpha)$ of the arclengths.

{\bf Proof.} ($\Rightarrow$) Let $\psi_\alpha(s_\alpha)$ and $\psi_\beta(s_\beta)$ be two regular similar curves with variable transformations in $\e^3$. Then there exists a variable transformation of the arclengths such that the tangent vectors and the principal normal vectors are the same (definition \ref{df-main1} and theorem \ref{th-03}). From equations (\ref{u3}) and (\ref{u4}) we have
\begin{equation}\label{u12}
\b_\beta(s_\beta)=\t_\beta(s_\beta)\times\n_\beta(s_\beta)=\t_\alpha(s_\alpha)\times\n_\alpha(s_\alpha)=\b_\alpha(s_\alpha).
\end{equation}

($\Leftarrow$) Let $\psi_\alpha(s_\alpha)$ and $\psi_\beta(s_\beta)$ be two regular curves in $\e^3$ which the same binormal vector under the arbitrary variable transformation $s_\beta=s_\beta(s_\alpha)$ of the arclengths. If we differentiate the equation (\ref{u11}) with respect to $s_\beta$ we have
\begin{equation}\label{u13}
-\tau_\beta(s_\beta)\n_\beta(s_\beta)=-\tau_\alpha(s_\alpha)\n_\alpha(s_\alpha)\frac{ds_\alpha}{ds_\beta}.
\end{equation}
The above equation leads to the following two equations
\begin{equation}\label{u14}
\n_\beta(s_\beta)=\n_\alpha(s_\alpha),
\end{equation}
From equations (\ref{u11}) and (\ref{u14}) we have
\begin{equation}\label{u16}
\t_\beta(s_\beta)=\n_\beta(s_\beta)\times\b_\beta(s_\beta)=\n_\alpha(s_\alpha)\times\b_\alpha(s_\alpha)=\t_\alpha(s_\alpha).
\end{equation}
The proof is complete.

\begin{theorem}\label{th05} Let $\psi_\alpha(s_\alpha)$ and $\psi_\beta(s_\beta)$ be two regular curves in $\e^3$. Then $\psi_\alpha(s_\alpha)$ and $\psi_\beta(s_\beta)$ are two similar curves with variable transformation if and only if the ratios of torsion and curvature are the same for all curves
\begin{equation}\label{u17}
\frac{\tau_\beta(s_\beta)}{\kappa_\beta(s_\beta)}=\frac{\tau_\alpha(s_\alpha)}{\kappa_\alpha(s_\alpha)},
\end{equation}
under the particular variable transformations ($\lambda^\beta_\alpha=\frac{ds_\beta}{ds_\alpha}=\frac{\kappa_\alpha}{\kappa_\beta}$) keeping equal total curvatures, i.e.,
\begin{equation}\label{u171}
\theta_\beta(s_\beta)=\int\,\kappa_\beta\,ds_\beta=\int\,\kappa_\alpha\,ds_\alpha=\theta_\alpha(s_\alpha).
\end{equation}
of the arclengths.
\end{theorem}

{\bf Proof.} Let $\psi_\alpha(s_\alpha)$ and $\psi_\alpha(s_\beta)$ are two similar curves with variable transformation. Then there exists a variable transformation of the arclengths such that the tangent and the binormal vectors are the same (definition \ref{df-main1} and theorem \ref{th04}). Differentiating the equations (\ref{u3}) and (\ref{u11}) we have
\begin{equation}\label{u18}
\kappa_\beta(s_\beta)\n_\beta(s_\beta)=\kappa_\alpha(s_\alpha)\n_\alpha(s_\alpha)\frac{ds_\alpha}{ds_\beta},
\end{equation}
\begin{equation}\label{u19}
-\tau_\beta(s_\beta)\n_\beta(s_\beta)=-\tau_\alpha(s_\alpha)\n_\alpha(s_\alpha)\frac{ds_\alpha}{ds_\beta}.
\end{equation}
which leads to the following two equations
\begin{equation}\label{u201}
\kappa_\beta(s_\beta)=\kappa_\alpha(s_\alpha)\frac{ds_\alpha}{ds_\beta}.
\end{equation}
\begin{equation}\label{u211}
\tau_\beta(s_\beta)=\tau_\alpha(s_\alpha)\frac{ds_\alpha}{ds_\beta}.
\end{equation}
The variable transformation (\ref{u17}) is the equation (\ref{u201}) after integration.  Dividing the above two equations (\ref{u211}) and (\ref{u201}), we obtain the equation (\ref{u17}) under the variable transformations (\ref{u171}).

($\Leftarrow$) Let $\psi_\alpha(s_\alpha)$ and $\psi_\beta(s_\beta)$ be two curves such that the equation (\ref{u17}) is satisfied under the variable transformation (\ref{u171}) of the arclengths. From theorem (\ref{th01}), the tangent vectors $\t_\alpha(s_\alpha)$ and $\t_\beta(s_\beta)$ of the two curves satisfy vector differential equations of third order as follows:
\begin{equation}\label{u22}
\Big(\frac{1}{f_\alpha(\theta_\alpha)}\t_\alpha''(\theta_\alpha)\Big)'+\Big(\frac{1+f_\alpha^2(\theta_\alpha)}{f(\theta_\alpha)}
\Big)\t_\alpha'(\theta_\alpha)
-\frac{f_\alpha'(\theta_\alpha)}{f_\alpha(\theta_\alpha)}\t_\alpha(\theta_\alpha)=0,
\end{equation}
\begin{equation}\label{u23}
\Big(\frac{1}{f_\beta(\theta_\beta)}\t_\beta''(\theta_\beta)\Big)'+\Big(\frac{1+f_\beta^2(\theta_\beta)}{f(\theta_\beta)}
\Big)\t_\beta'(\theta_\beta)
-\frac{f_\beta'(\theta_\beta)}{f_\beta(\theta_\beta)}\t_\beta(\theta_\beta)=0,
\end{equation}
where $f_\alpha(\theta_\alpha)=\frac{\tau_\alpha(\theta_\alpha)}{\kappa_\alpha(\theta_\alpha)}$, $f_\beta(\theta_\beta)=\frac{\tau_\beta(\theta_\beta)}{\kappa_\beta(\theta_\beta)}$,  $\theta_\alpha=\int\,\kappa_\alpha(s_\alpha)\,ds_\alpha$ and $\theta_\beta=\int\,\kappa_\beta(s_\beta)\,ds_\beta$.

The equation (\ref{u17}) leads to
\begin{equation}\label{u24}
f_\beta(\theta_\beta)=f_\alpha(\theta_\alpha),
\end{equation}
under the variable transformations $\theta_\beta=\theta_\alpha$. So that the two equations (\ref{u22}) and (\ref{u23}) under the equation (\ref{u17}) and the transformation (\ref{u171}) are the same. Hence the solution is the same, i.e., the tangent vectors are the same which completes the proof of the theorem.

\section{New classifications of curves}

In this section, we will apply our definition of similar curves with variable transformations to deduce the position vectors of some special curves. First we can call the two curves $\psi_\alpha(s_\alpha)$ and $\psi_\beta(s_\beta)$  similar curves with variable transformation if and only if there exists an arbitrary function $\lambda_\beta^\alpha=\frac{ds_\alpha}{ds_\beta}$ such that the curvature and torsion of the curve $\psi_\beta$ are the curvature and torsion of the curve $\psi_\alpha$ multiplied by this arbitrary function i,e.,
\begin{equation}\label{u25}
\kappa_\beta=\kappa_\alpha\,\lambda_\beta^\alpha,\,\,\,\,\,
\tau_\beta=\tau_\alpha\,\lambda_\beta^\alpha.
\end{equation}

{\bf Class 1.} If the curve is straight line then the curvature is $\kappa=0$. Under the variable transformation $\lambda$ the curvature does not change. So we have the following lemma:

\begin{lemma}\label{lm-1}
The straight line alone forms a family of similar curves with variable transformation.
\end{lemma}

{\bf Class 2.} If the curve is a plane curve then the torsion is $\tau_\alpha=0$. Under the variable transformation $\lambda$ the torsion does not change. So we have the following lemma:

\begin{lemma}\label{lm-2}
The family of plane curves forms a family of similar curves with variable transformations.
\end{lemma}

We can deduce the position vector of a plane curve using the definition of similar curves with variable transformation as follows:

The simplest example of a plane curve is a circle of radius 1. The natural representation of this circle can be written in the form:
\begin{equation}\label{u30}
\psi_\alpha(u)=\Big(\sin[u],-\cos[u],0\Big),
\end{equation}
where $s_\alpha=u$ is the arclength of the circle and the curvature is $\kappa_\alpha(u)=1$. The tangent vector of this circle takes the form:
\begin{equation}\label{u31}
\t_\alpha(u)=\Big(\cos[u],\sin[u],0\Big).
\end{equation}
From theorem (\ref{th02}) we can write any plane curve as the following:
\begin{equation}\label{u32}
\psi_\beta(s)=\int\Big(\cos\Big[u[s]\Big],\sin\Big[u[s]\Big],0\Big)ds.
\end{equation}
where $s_\beta=s$. From the equation (\ref{u25}), we have
\begin{equation}\label{u33}
ds_\alpha=\lambda^\alpha_\beta\,ds_\beta=\frac{\kappa_\beta}{\kappa_\alpha}\,ds_\beta.
\end{equation}
or
\begin{equation}\label{u34}
s_\alpha(s_\beta)=\int\,\frac{\kappa_\alpha}{\kappa_\beta}\,ds_\beta.
\end{equation}
If we put the curvature $\kappa_\beta=\kappa(s)$ $(s_\beta=s)$, we have
\begin{equation}\label{u35}
u(s)=\int\,\kappa(s)\,ds.
\end{equation}
Then the position vector of the plane curve with arbitrary curvature $\kappa(s)$ takes the following form:
\begin{equation}\label{u36}
\psi(s)=\int\Big(\cos\Big[\int\,\kappa(s)\,ds\Big],\sin\Big[\int\,\kappa(s)\,ds\Big],0\Big)ds.
\end{equation}
which is the position vector of a plane curve introduced in \cite{lips}.

{\bf Class 3.} If the curve $\psi_\alpha$ is a general helix $(\frac{\tau_\alpha}{\kappa_\alpha}=m)$, where $m$ is a constant, in the form $m=\cot[\phi]$ and $\phi$ is the angle between the tangent vector and the axis of the helix. Then any similar curve $\psi_\beta$ with this helix has the the property  $\frac{\tau_\beta}{\kappa_\beta}=m$. So that we have the following lemma:

\begin{lemma}\label{lm-3}
The family of general helices with fixed angle $\phi$ between the axis of a general helix and the tangent vector forms a family of similar curves with variable transformations.
\end{lemma}

We can deduce the position vector of a general helix using the definition of similar curves with variable transformations as follows:

The simplest example of a general helix is a circular helix or $W$-curve. The natural representation of a circular helix is:
\begin{equation}\label{u40}
\psi_\alpha(u)=\Big(\sqrt{1-n^2}\,\sin[u],-\sqrt{1-n^2}\,\cos[u],n\,u\Big),
\end{equation}
where $u$ is the arclength of the circular helix and $n=\cos[\phi]$, where $\phi$ is the constant angle between the tangent vector and the axis of a circular helix. The curvature of this circular helix is $\kappa_\alpha(u)=\sqrt{1-n^2}$. The tangent vector of this curve takes the form:
\begin{equation}\label{u41}
\t_\alpha(u)=\Big(\sqrt{1-n^2}\cos[u],\sqrt{1-n^2}\sin[u],n\Big).
\end{equation}
From theorem (\ref{th02}) we can write any general helix as the following:
\begin{equation}\label{u42}
\psi_\beta(s)=\int\Big(\sqrt{1-n^2}\cos\Big[u(s)\Big],\sqrt{1-n^2}\sin\Big[u(s)\Big],n\Big)ds.
\end{equation}
From equation (\ref{u34}) we have
\begin{equation}\label{u44}
u(s)=\int\,\frac{\kappa(s)}{\sqrt{1-n^2}}\,ds.
\end{equation}
where $\kappa_\beta=\kappa(s)$, $(s_\beta=s)$. Then the position vector of the general helix with arbitrary curvature $\kappa(s)$ takes the following form:
\begin{equation}\label{u45}
\psi=\int\Big(\sqrt{1-n^2}\cos\Big[\int\frac{\kappa(s)}{\sqrt{1-n^2}}ds\Big],
\sqrt{1-n^2}\sin\Big[\int\frac{\kappa(s)}{\sqrt{1-n^2}}ds\Big],n\Big)ds.
\end{equation}
which is the position vector of a general helix introduced in \cite{ali2}.

{\bf Class 4.} If the curve is a slant helix then the relation (\ref{u02}) between the torsion and curvature is satisfied. Let $\psi_\alpha$ and $\psi_\beta$ be two slant helices such that the transformation (\ref{u171}) is satisfied. Using the relation (\ref{u02}) and (\ref{u171}), it is easy to prove that:
$$
\frac{\tau_\beta}{\kappa_\beta}=\frac{m\,\theta_\beta}{\sqrt{1-m^2\theta^2_\beta}}=
\frac{m\,\theta_\alpha}{\sqrt{1-m^2\theta^2_\alpha}}=\frac{\tau_\alpha}{\kappa_\alpha},
$$
where $m$ is a constant value, $m=\cot[\phi]$ and $\phi$ is the angle between the principal normal vector and the axis of a slant helix. So that we have the following lemma:

\begin{lemma}\label{lm-4}
The family of a slant helices with fixed angle $\phi$ between the axis of a slant helix and the principal normal vector forms a family of similar curves with variable transformation.
\end{lemma}

Now, we can deduce the position vector of a slant helix using the definition of similar curves with variable transformations as follows:

The simplest example of a slant helix is Salkowski curve \cite{ali3, mont1, salkow}. The explicit parametric  representation of a Salkowski curve $\psi_\alpha(u)=\big(\psi_1(u),\psi_2(u),\psi_3(u)\big)$ takes the form:
\begin{equation}\label{u50}
\left\{
\begin{array}{ll}
\psi_1(u)=-\frac{n}{4m}\Big[\frac{n-1}{2n+1}\cos[(2n+1)t]+
\frac{n+1}{2n-1}\cos[(2n-1)t]-2\cos[t]\Big],\\
\psi_2(u)=-\frac{n}{4m}\Big[\frac{n-1}{2n+1}\sin[(2n+1)t]-
\frac{n+1}{2n-1}\sin[(2n-1)t]-2\sin[t]\Big],\\
\psi_3(u)=\frac{n}{4m^2}\cos[2nt],
\end{array}
\right.
\end{equation}
where $t=\frac{1}{n}\mathrm{arcsin}(mu)$, $m=\frac{n}{\sqrt{1-n^2}}$, $n=\cos[\phi]$ and $\phi$ is the constant angle between the axis of a slant helix and the principal normal vector. The curvature of the above curve is $1$ and the torsion is
$$
\tau(u)=\tan[nt]=\frac{m\,u}{\sqrt{1-m^2\,u^2}}.
$$
It is worth noting that: the variable $t$ is a parameter while the variable $u$ is the natural parameter.

The tangent and the principal normal vectors of the Salkowski curve (\ref{u50}) take the forms:
\begin{equation}\label{u51}
\t_\alpha(u)=-\Big(n\cos[t]\sin[nt]-\sin[t]\cos[nt], n\sin[t]\sin[nt]+\cos[t]\cos[nt],\frac{n}{m}\sin[nt]\Big).
\end{equation}
\begin{equation}\label{u52}
\n_\alpha(u)=\Big(\sqrt{1-n^2}\cos[t],\sqrt{1-n^2}\sin[t],n\Big),
\end{equation}
It is easy to write the tangent vector (\ref{u51}) in the simple form:
\begin{equation}\label{u53}
\t_\alpha(u)=\int\n(u)du=\int\Big(\sqrt{1-n^2}\cos[t],\sqrt{1-n^2}\sin[t],n\Big)du.
\end{equation}
From theorem (\ref{th02}) we can write any slant helix as the following:
\begin{equation}\label{u54}
\psi_\beta(s)=\int\Big[\int\Big(\sqrt{1-n^2}\cos[t],\sqrt{1-n^2}\sin[t],n\Big)du\Big]ds.
\end{equation}
From equation (\ref{u34}) we have
\begin{equation}\label{u55}
u(s)=\int\,\kappa(s)\,ds,\,\,\,\,\,du=\kappa(s)\,ds,
\end{equation}
where $\kappa_\beta=\kappa(s)$, $(s_\beta=s)$. Substituting equation (\ref{u55}) in (\ref{u54}) we obtain the position vector of a similar curve $\psi_\beta(s)=\big(\psi_1(s),\psi_2(s),\psi_3(s)\big)$of a slant helix with arbitrary curvature $\kappa(s)$ as follows:
\begin{equation}\label{u56}
\left\{
\begin{array}{ll}
\psi_1(s)=\frac{n}{m}\,\int\Bigg[\int\kappa(s)\cos\Big[\frac{1}{n}
          \mathrm{arcsin}\Big(m\int\kappa(s)ds\Big)\Big]ds\Bigg]ds,\\
\psi_2(s)=\frac{n}{m}\,\int\Bigg[\int\kappa(s)\sin\Big[\frac{1}{n}
          \mathrm{arcsin}\Big(m\int\kappa(s)ds\Big)\Big]ds\Bigg]ds,\\
\psi_3(s)=n\,\int\Big[\int\kappa(s)ds\Big]ds,
\end{array}
\right.
\end{equation}
which is the position vector of a slant helix introduced in \cite{ali3}.

We hope that we can introduce new classes of similar curves and deduce the position vectors of these classes in future work.


\end{document}